\documentclass[11pt]{article}
\usepackage{amsfonts,amssymb,color}
\usepackage[mathscr]{eucal}
\usepackage{amsmath, amsthm}
\usepackage{mathrsfs}
\usepackage{amsbsy}
\usepackage{dsfont}
\usepackage{bbm}
\usepackage{wasysym}
\usepackage{stmaryrd}
\usepackage{upgreek}
\usepackage{hyperref}
\input xypic
\xyoption{all}
\begin{document}
\baselineskip = 5mm
\newcommand \ZZ {{\mathbb Z}}
\newcommand \FF {{\mathbb F}}
\newcommand \NN {{\mathbb N}}
\newcommand \QQ {{\mathbb Q}}
\newcommand \RR {{\mathbb R}}
\newcommand \CC {{\mathbb C}}
\newcommand \PR {{\mathbb P}}
\newcommand \AF {{\mathbb A}}
\newcommand \ud {{\Bbbk }}
\newcommand \bgf {{\mathbb K}}
\newcommand \bcA {{\mathscr A}}
\newcommand \bcB {{\mathscr B}}
\newcommand \bcC {{\mathscr C}}
\newcommand \bcD {{\mathscr D}}
\newcommand \bcE {{\mathscr E}}
\newcommand \bcF {{\mathscr F}}
\newcommand \bcG {{\mathscr G}}
\newcommand \bcH {{\mathscr H}}
\newcommand \bcM {{\mathscr M}}
\newcommand \bcN {{\mathscr N}}
\newcommand \bcI {{\mathscr I}}
\newcommand \bcJ {{\mathscr J}}
\newcommand \bcK {{\mathscr K}}
\newcommand \bcL {{\mathscr L}}
\newcommand \bcO {{\mathscr O}}
\newcommand \bcP {{\mathscr P}}
\newcommand \bcQ {{\mathscr Q}}
\newcommand \bcR {{\mathscr R}}
\newcommand \bcS {{\mathscr S}}
\newcommand \bcT {{\mathscr T}}
\newcommand \bcU {{\mathscr U}}
\newcommand \bcV {{\mathscr V}}
\newcommand \bcW {{\mathscr W}}
\newcommand \bcX {{\mathscr X}}
\newcommand \bcY {{\mathscr Y}}
\newcommand \bcZ {{\mathscr Z}}
\newcommand \ga {{\mathfrak a}}
\newcommand \gb {{\mathfrak b}}
\newcommand \gc {{\mathfrak c}}
\newcommand \gd {{\mathfrak d}}
\newcommand \gm {{\mathfrak m}}
\newcommand \gn {{\mathfrak n}}
\newcommand \gp {{\mathfrak p}}
\newcommand \gq {{\mathfrak q}}
\newcommand \gQ {{\mathfrak Q}}
\newcommand \gP {{\mathfrak P}}
\newcommand \gT {{\mathfrak T}}
\newcommand \gC {{\mathfrak C}}
\newcommand \gD {{\mathfrak D}}
\newcommand \gM {{\mathfrak M}}
\newcommand \gS {{\mathfrak S}}
\newcommand \gH {{\mathfrak H}}
\newcommand \gA {{\rm {A}}}
\newcommand \gB {{\rm {B}}}
\newcommand \catC {{\sf C}}
\newcommand \catD {{\sf D}}
\newcommand \Sets {{\sf Sets}}
\newcommand \Sch {{\sf Sch }}
\newcommand \Vect {{\sf Vect}}
\newcommand \SmProj {{\sf SmProj}}
\newcommand \MM {{\sf MM}}
\newcommand \HS {{\sf HS}}
\newcommand \MHS {{\sf MHS}}
\newcommand \Zar {\rm {Zar}}
\newcommand \Nis {\rm {Nis}}
\newcommand \cdh {\rm {cdh}}
\newcommand \h {\rm {h}}
\newcommand \fppf {{\rm {fppf}}}
\newcommand \et {\rm {\acute e t}}
\newcommand \Mot {{\sf Mot}}
\newcommand \CHM {{\sf Chow}}
\newcommand \Num {{\sf Num}}
\newcommand \DM {{\sf DM}}
\newcommand \Ta {{\mathbb T}}
\newcommand \uno {{\mathbbm {1}}}
\newcommand \Le {{\mathbb L}}
\newcommand \ptr {{\pi _2^{\rm tr}}}
\newcommand \Spec {{\rm {Spec}}}
\newcommand \bSpec {{\bf {Spec}}}
\newcommand \Pic {{\rm {Pic}}}
\newcommand \PicF {{\it {Pic}}}
\newcommand \PicS {{\rm {Pic}}}
\newcommand \Jac {{{J}}}
\newcommand \AlbS {{\rm {Alb}}}
\newcommand \alb {{\rm {alb}}}
\newcommand \NS {{{NS}}}
\newcommand \Corr {{Corr}}
\newcommand \Sym {{\rm {Sym}}} 
\newcommand \Symcl {{\rm {S}}} 
\newcommand \Alt {{\rm {Alt}}}
\newcommand \Prym {{\rm {Prym}}}
\newcommand \HilbF {{\it Hilb}}
\newcommand \HilbS {{\it Hilb}}
\newcommand \Mor {{\rm Mor}}
\newcommand \MorF {{\it Mor}}
\newcommand \MorS {{\it Mor}}
\newcommand \pdiv {{\rm {div}}}
\newcommand \Bl {{\rm {Bl}}}
\newcommand \codim {{\rm {codim}}}
\newcommand \Proj {{\rm {Proj}}}
\newcommand \bProj {{\bf {Proj}}}
\newcommand \Div {{\rm {Div}}}
\newcommand \prim {{\rm {prim}}}
\newcommand \stalk {{\rm st}}
\newcommand \seminorm {{\rm {sn}}}
\newcommand \hc {{\rm hc}} 
\newcommand \BC {{\rm {BC}}}
\newcommand \cnj {{\rm {C}}}
\newcommand \Res {{\rm {Res}}}
\newcommand \card {{\rm {card}}}
\newcommand \cone {{\rm {cone}}}
\newcommand \cha {{\rm {char}}}
\newcommand \eff {{\rm {eff}}}
\newcommand \cl {{\rm {cl}}}
\newcommand \tr {{\rm {tr}}}
\newcommand \pr {{\rm {pr}}}
\newcommand \ev {{\rm {ev}}}
\newcommand \sep {{\rm {sep}}}
\newcommand \alg {{\rm {alg}}}
\newcommand \im {{\rm im}}
\newcommand \rd {{\rm {red}}}
\newcommand \op {{\rm op}}
\newcommand \Hom {{\rm Hom}}
\newcommand \uHom {{\underline {\rm Hom}}}
\newcommand \cHom {{\mathscr H\! }{\it om}}
\newcommand \Ext {{\rm Ext}}
\newcommand \cExt {{\mathscr E\! }{\it xt}}
\newcommand \colim {{{\rm colim}\, }} 
\newcommand \End {{\rm {End}}}
\newcommand \coker {{\rm {coker}}}
\newcommand \id {{\rm {id}}}
\newcommand \van {{\rm {van}}}
\newcommand \spc {{\rm {sp}}}
\newcommand \Ob {{\rm Ob}}
\newcommand \Aut {{\rm Aut}}
\newcommand \cor {{\rm {cor}}}
\newcommand \res {{\rm {res}}}
\newcommand \Gal {{\rm {Gal}}}
\newcommand \PGL {{\rm {PGL}}}
\newcommand \Gr {{\rm {Gr}}}
\newcommand \Tor {{\rm {Tor}}}
\newcommand \Sing {{\rm {Sing}}}
\newcommand \spn {{\rm {span}}}
\newcommand \univ {{\rm {\, univ}}}
\newcommand \Nm {{\rm {Norm}}}
\newcommand \inv {{\rm {inv}}}
\newcommand \even {{\rm {even}}}
\newcommand \md {{\rm {mod}\, }}
\newcommand \ind {{\rm {ind}}}
\newcommand \Gm {{{\mathbb G}_{\rm m}}}
\newcommand \trdeg {{\rm {tr.deg}}}
\newcommand \sing {{\rm {sing}}}
\newcommand \tame {\rm {tame }}
\newcommand \eq {{\rm {eq}}}
\newcommand \length {{\rm {length}}}
\newcommand \ord {{\rm {ord}}}
\newcommand \equi {{\rm {equi}}}
\newcommand \ab {{\rm {ab}}}
\newcommand \Fix {{\rm {Fix}}}
\newcommand \trp {{\rm {t}}}
\newcommand \cat {{\rm {cat}}}
\newcommand \aff {{\rm {aff}}}
\newcommand \Const {{\rm {Const}}}
\newcommand \num {{\rm {num}}}
\newcommand \rat {{\rm rat}}
\newcommand \tors {{\rm {tors}}}
\newcommand \coeq {{{\rm coeq}\, }}
\newcommand \supp {{\rm Supp}}
\newcommand \sm {{\rm sm}}
\newcommand \reg {{\rm reg}}
\newcommand \var {{\rm var}}
\newcommand \norm {{\rm {nor}}}
\newcommand \Gys {{\rm {Gys}}}
\newcommand \Fr {{\rm {Fr}}}
\newcommand \hol {{\rm {h}}}
\newcommand \bideg {{\rm {bideg}}}
\newcommand \loc {{\rm {loc}}}
\newcommand \znak {{\natural }}
\newcommand \znakk {{\sharp }}
\newcommand \znakkk {{\flat }}
\newcommand \qand {{\quad \hbox{and}\quad }}
\newcommand \qqand {{\quad \quad \hbox{and}\quad \quad }}
\newcommand \qqqand {{\quad \quad \quad \hbox{and}\quad \quad \quad }}
\newcommand \heither {{\hbox{either}\quad }}
\newcommand \qor {{\quad \hbox{or}\quad }}
\newcommand \qqor {{\quad \quad \hbox{or}\quad \quad }}
\newcommand \lra {\longrightarrow}
\newcommand \hra {\hookrightarrow}
\def\blue {\color{blue}}
\def\red {\color{red}}
\def\green {\color{green}}
\newtheorem{theorem}[subsubsection]{Theorem}
\newtheorem{lemma}[subsubsection]{Lemma}
\newtheorem{corollary}[subsubsection]{Corollary}
\newtheorem{proposition}[subsubsection]{Proposition}
\newtheorem{remark}[subsubsection]{Remark}
\newtheorem{definition}[subsubsection]{Definition}
\newtheorem{conjecture}[subsubsection]{Conjecture}
\newtheorem{example}[subsubsection]{Example}
\newtheorem{question}[subsubsection]{Question}
\newtheorem{comment}[subsubsection]{Comment}
\newtheorem{assumption}[subsubsection]{Assumption}
\newtheorem{fact}[subsubsection]{Fact}
\newtheorem{crucialquestion}[subsubsection]{Crucial Question}
\newtheorem{claim}[subsubsection]{Claim}
\newtheorem{terminology}[subsubsection]{Terminology}
\newtheorem{aim}[subsubsection]{Aim}
\newenvironment{pf}{\par\noindent{\em Proof}.}{\hfill\framebox(6,6)
\par\medskip}
\title{\bf Motivic obstruction to rationality \\ of a general cubic in $\PR ^5$}



\author{\sc Vladimir Guletski\u \i }


\maketitle

\begin{abstract}
\noindent We introduce integrally essentially indecomposable motives and prove that, if the integral motive of a smooth projective surface over a field of sufficiently big transcendence degree is integrally essentially indecomposable, then a very general cubic fourfold in $\PR ^5$ is not rational. We also prove a lifting theorem saying that, given a smooth projective family of surfaces over a Henselian DVR, if the motive of the special fibre is integrally essentially indecomposable, then so is the motive of the generic fibre. This suggests a possible reduction of the cubic fourfold conjecture to certain arithmetic phenomena in positive characteristic.
\end{abstract}
















\tableofcontents

\section{Introduction}
\label{intro}

A long standing conjecture in algebraic geometry says that a very general cubic hypersurface in $\PR ^5$ is not rational. Since such fourfolds are unirational, the conjecture is a particular case of the L\"uroth problem. Whereas the L\"uroth problem for cubic threefolds was solved in terms of abelian invariants (Prym varieties, theta-divisors, etc), \cite{ClemensGriffiths}, the numerous attempts to develop an analog of the Clemens-Griffiths theory, which could be appropriate in dimension $4$, have not achieved the desired result yet (but see the recent progress in \cite{KKPYY}). 

The aim of this paper is to develop a motivic version of Kulikov's reduction in \cite{Kulikov}, given in terms of integral Hodge structures, and think of it in arithmetic perspective.

Recall that a natural birational invariant of cycle-theoretic nature is the Chow group of $0$-cycles modulo rational equivalence on a variety over a field. Important developments along this line include the notion of $CH_0$-triviality introduced in \cite{ACP}. In \cite{VoisinUnirat} Voisin proved that $CH_0$-nontriviality is a deformable property in families, and used this to prove the stable non-rationality for the desingularization of a very general quartic double solid with at most seven nodes. In \cite{Colliot-ThelenePirutka} Colliot-Th\'el\`ene and Pirutka used similar method to prove the existence of not stably rational smooth quartic hypersurfaces in $\PR ^4$. However, as we do not know a single example of a nonrational cubic fourfold in $\PR ^5$, it is not clear how to use the deformation of $CH_0$-nontriviality in the striking dimension $4$ case. The goal of our paper is to convince the reader that the cycle-theoretic reduction should be to certain phenomena of arithmetic nature in positive characteristic. 

 Kulikov's reduction says that if the integral trenscendental Hodge structure of a smooth projective surface over $\CC $ is indecomposable, then the cubic fourfold conjecture follows, see \cite{Kulikov}. Later on, Auel, B\"ohning and Graf von Bothmer see \cite{ABG} demonstrated, using Macalauy2 computation, that the transcendental Hodge structure of the Fermat sextic surface in $\PR ^3$ is decomposable, which seemingly dismiss Kulikov's approach. Our point is that it does not. 

Indeed, the minimal field of definition of the surface $S$ whose integral transcendental Hodge structure decomposes in Lemma 2 on page 61 in \cite{Kulikov} is not the field $\QQ $, as we have have it in case of Fermat sextic, but rather $S$ has a model over a field of positive transcendence degree over the primary subfield. It means that if we add the requirement saying that $S$ is sufficiently general, or, equivalently, its minimal field of definition, the sense of Weil, \cite{Weil}, is big enough, in the Nondecomposability Conjecture on page 59 in \cite{Kulikov}, then Kulikov's reduction is still valid. 

But we can do more. Instead of looking at the integral transcendental Hodge structures of a smooth projective surface $S$ over $\CC $, one can work with the transcendental motive $M^2_{\tr }(S)$ of the surface $S$. Of course, it is defined rationally, but the analysis of this motive in \cite{KMP} reveals that any reasonable concept of integral indecomposability of $M^2_{\tr }(S)$ must be equivalent to what we call essential integral indecomposability of the whole motive $M(X)$. 

This all leads us to a natural motivic analog of Kulikov's Hodge-theoretic obstruction to rationality of a very general cubic fourfold in $\PR ^5$, which would avoid the difficulties above. Moreover, here are two advantages of the motivic approach advocated in our paper. The first one is that there is no phantom submotives in a motive, provided it is finite-dimensional, see Proposition 7.5 in \cite{Kimura}. The second advantage is that the obstruction to rationality of a fourfold is given in terms of rational equivalence of $0$-cycles on surfaces, rather than on the fourfold itself, and this allows us to reduce the cubic fourfold problem to positive characteristic. 


\begin{itemize}

\item[]{}

{\sc Motivic indecomposability conjecture.} {\it Let $k$ be a field of characteristic $0$, let $S$ be a smooth projective surface over $k$, and let $k_0$ is the minimal field of definition of the surface $S$. If 
  $$
  \trdeg (k_0/\QQ )\gg 0\; ,
  $$
the motive $M(S)$ is integrally essentially indecomposable.}

\end{itemize}

\medskip

Now recall also the well-know conjecture due to Kimura and O'Sullivan which asserts that all Chow motives are finite-dimensional, see \cite{Andre} and \cite{AndrePanorama}. In the paper we prove the following conditional

\medskip

\begin{itemize}

\item[]{}

{\sc Theorem A.} {\it If the motivic indecomposability conjecture is true, and if the motive of any smooth projective surface is finite-dimensional, then a very general cubic fourfold hypersurface in $\PR ^5$ is not rational.}

\end{itemize}

\medskip

Our second result is the lifting theorem. Let $B$ be the spectrum of a Henselian discrete valuation ring, let $0\in B$ be the closed point of $B$, and let $\eta $ be its generic point. Let also $X/B$ be a smooth projective family of surfaces over $B$, let $X_0$ be the closed fibre, and let $X_{\eta }$ be the generic one.

\medskip

\begin{itemize}

\item[]{}

{\sc Theorem B.} {\it If the motive $M(X_0)$ is integrally essentially indecomposable, so is the motive $M(X_{\eta })$.}

\end{itemize}

\medskip

We believe that integrally essentially indecomposable motives are 'motivic atoms' whose enriched Hodge realizations are the Hodge atoms in \cite{KKPYY}. 

In the last section we discuss what kind of arithmetic phenomena can cause integral essential indecomposability of motives in characteristic $p>0$, and how to use it  in order to approach the motivic indecomposability conjecture above. 

\bigskip

{\sc Acknowledgements.} I am grateful to Alexander Kuznetsov, Mingmin Shen and Alexander Tikhomirov for useful discussions. I am also grateful to the Center for Geometry and Physics at the Institute for Basic Science in Pohang (South Korea), where the first version of this paper was written long ago, and to the organizers of the VII Escuela IMCA conference in Algebraic Geometry (September 2019) where I reported Theorem B for the first time.

\medskip 

\section{Integral motives}

\subsection{Basic definitions}
\label{basics}

Let $k$ be an algebraically closed field. If $X$ is a scheme of finite type over $k$, let $CH_r(X)$ be the Chow group of dimension $r$ algebraic cycles modulo rational equivalence on $X$ with coefficients in $\ZZ $, and let also $A_r(X)$ be the subgroup generated by algebraically trivial cycle classes in $CH_r(X)$. If $X$ is equidimensional of dimension $n$, then $CH^{n-r}(X)$ and $A^{n-r}(X)$ stay for $CH_r(X)$ and $A_r(X)$ respectively. For any abelian group $A$ we denote by $A_{\QQ }$ the tensor product $A\otimes _{\ZZ }\QQ $.

The category $\CHM (k,\ZZ )$ of Chow motives over $k$ with coefficients in $\ZZ $ can be defined in a usual way. If $X$ and $Y$ are two smooth projective varieties over $k$, and 
  $$
  X=\cup _jX_j
  $$ 
is the decomposition of $X$ into connected components, then 
  $$
  CH^m(X,Y)=\oplus CH^{n_j+m}(X_j\times Y)
  $$
is the group of correspondences of degree $m$ from $X$ to $Y$, where $n_j$ is the dimension of the component $X_j$. For any two correspondences 
  $$
  \alpha \in CH^m(X,Y)\qqand \beta \in CH^n(Y,Z)\; ,
  $$ 
 their composition is given by the standard pullback-intersection-pushforward formula
   $$
   \beta \circ \alpha ={p_{13}}_*(p_{12}^*(\alpha )\cdot p_{23}^*(\beta ))
   $$
 where the central dot stays for the intersection of cycle classes and the projections from $X\times Y\times Z$ onto $X\times Y$ and $Y\times Z$ are obvious. The composition is bilinear, and we obtain the homomorphism 
   $$
   CH^n(Y,Z)\otimes CH^m(X,Y)\to CH^{m+n}(X,Z)\; .
   $$

The objects of $\CHM (k)$ are triples 
  $$
  (X,\Sigma ,m)\; ,
  $$
where $X$ is a smooth projective variety over $k$,
  $$
  \Sigma \in CH^0(X,X)
  $$ 
is an idempotent, sometimes called projector, in the associative ring $CH^0(X,X)$ of degree $0$ correspondences from $X$ to $X$, and $m$ is an integer. 

For two motives 
  $$
  M=(X,\Sigma ,m)\qqand N=(Y,\Xi ,n)\; ,
  $$
the group 
  $$
  \Hom _{\CHM (k)}(M,N)=\{ \Xi \circ \Phi \circ \Sigma \; |\; \Phi \in CH^{n-m}(X,Y)\} 
  $$ 
 consists of all triple compositions $\Xi \circ \Phi \circ \Sigma $, where $\Phi $ is of the right degree $n-m$. 
 
 Let $\SmProj (k)$ be the category of smooth projective varieties over $k$. The contravariant functor 
   $$
   M:\SmProj (k)\to \CHM (k,\ZZ )
   $$
 sends every smooth projective variety $X$ over $k$ to its motive
   $$
   M(X)=(X,\Delta _X,0)\; ,
   $$
 where $\Delta _X$ is the diagonal class in $CH^0(X,X)$. If
   $$
   f:X\to Y
   $$
 is a morphism of smooth projective varieties over $k$, then the transposed graph $\Gamma ^{\trp }_f$ of the morphism $f$ is an element of the group $CH^0(Y,X)$, and the graph of the identity map for $X$ to $X$ is the diagonal class $\Delta _X$, of course. Then
   $$
   M(f)=\Gamma ^{\trp }_f:M(Y)\to M(X)
   $$
 is the value of the functor $M$ on $f$. 
 
 If $\Sigma $ is an idempotent in $CH^0(X,X)$, it is sometime more convenient to write 
   $$
   M_{\Sigma }
   $$ 
instead of the triple $(X,\Sigma ,0)$.

It is important that the category $\CHM (k,\ZZ )$ is symmetric monoidal, where the monoidal product is given by the formula
  $$
  (X,\Sigma ,m)\otimes (Y,\Xi ,n)=(X\times Y,\Sigma \otimes \Xi ,m+n)\; .
  $$
The triple 
  $$
  \uno =(\Spec (k),\Delta ,0)
  $$ 
is the monoidal unit, and the triple 
  $$
  \Le =(\Spec (k),\Delta ,-1)
  $$ 
is called the Lefschetz motive over $k$. Then
  $$
  M(\PR ^1)=\uno \oplus \Le \; .
  $$
We will be also using the Tate motive 
  $$
  \Ta =\Le ^{-1}=(\Spec (k),\Delta ,1)\; ,
  $$
which is the monoidal inverse to the motive $\Le $ in $\CHM (k)$.

The category $\CHM (k,\QQ )$ with coefficients in $\QQ $ can be constructed in the same way, but using the $\QQ $-localized Chow groups 
  $$
  CH^n(X)_{\QQ }\; .
  $$
It is tensor rigid, in the sense of \cite{DeligneMilne}.

Similarly, one can construct the category 
  $$
  \Mot _{\sim }(k,R)
  $$
of pure motives over $k$ with coefficients in a commutative ring $R$ modulo any adequate equivalence relation $\sim $ on algebraic cycles, see \cite{Samuel} for the notion of an adequate equivalence relation. Then
  $$
  \CHM (k,\ZZ )=\Mot _{\rat }(k,\ZZ )\; ,
  $$
and the same for pure motives with coefficients in $\QQ $.

If 
  $$
  \Num (k,\ZZ )=\Mot _{\num }(k,\ZZ )
  $$ 
is the category of pure motives modulo numerical equivalence over $k$, we obtain the monoidal functor 
  $$
  \CHM (k,\ZZ )\to \Num (k,\ZZ )
  $$
sending a Chow motive $M=(X,\Sigma ,m)$ to 
  $$
  \bar M=(X,\bar \Sigma ,m)\; ,
  $$
where $\bar \Sigma $ is $\Sigma $ modulo numerical equivalence on $X$, and similarly on morphisms.

It is important that the category $\Num (k,\QQ )$ is semisimple abelian by Jannsen's result, see \cite{JannsenSemisimple}.

Now, in \cite{BlochSrinivas} Bloch and Srinivas studied important consequences of a certain assumption saying that a multiple of the diagonal $\Delta _X$ is rationally equivalent, with coefficients in $\ZZ $, to a linear combination of prime $n$-cycles on $X\times X$, each of which is not dominant on the left or right. Following \cite{Barbieri-Viale}, we say that an integral closed subscheme $Z$ of codimension $n$ in $X\times X$, i.e. a prime $n$-cycle, is balanced on the left (respectively, on the right) if there exists an equi-dimensional Zariski closed subscheme $Y\subset X$ with $\dim (Y)<n$, such that $Z$ is contained in $Y\times X$ (respectively, in $X\times Y$). In other words, $Z$ is balanced if it is not dominant on the left or on the right. 

A codimension $n$ cycle is balanced on $X\times X$ if it is rationally equivalent, with coefficients in $\ZZ $, to a sum of prime cycles, each of which is balanced on the left or right. A class in $CH^0(X,X)=CH^n(X\times X)$ is balanced if it can be represented by a balanced cycle on $X\times X$. 

Certainly, the same definitions can be also given with coefficients in $\QQ $, and modulo any adequate equivalence relation on algebraic cycles. We will use the corresponding terminology below with no additional definitions.  

It is essential for what follows that the subspace generated by balanced cycles is a two-sided ideal 
  $$
  BCH^n(X\times X)\subset CH^n(X\times X)
  $$
in the associative ring of degree $0$ correspondences from $X$ to $X$ with coefficients in $\ZZ $, see p. 309 in \cite{Fulton}. Accordingly, the $\QQ $-vector space $BCH^n(X\times X)_{\QQ }$ is an ideal in the associative algebra $CH^n(X\times X)_{\QQ }$. We may also write $BCH^0(X,X)$ and $BCH^0(X,X)_{\QQ }$ sometimes. 

Now, let $M=(X,\Sigma ,m)$ be a motive in $\CHM (k,\ZZ )$. We will say that $M$ decomposes essentially in $\CHM (k,\ZZ )$, if 
  $$
  \Sigma =\Lambda +\Xi 
  $$ 
is a sum of two mutually orthogonal projectors in $CH^0(X,X)$, each of which is not balanced modulo rational equivalence with coefficients in $\QQ $ on $X\times X$ and does not vanish modulo homological equivalence on $X\times X$, also with coefficients in $\QQ $. 

Otherwise, we will say that $M$ is integrally essentially indecomposable in $\CHM (k,\ZZ )$. The reader may wish to call integrally essentially indecomposable motives motivic atoms in $\CHM (k,\ZZ )$.

\subsection{The Chow-K\"unneth decomposition}

Choose and fix a prime $l$, different from the characteristic of the ground field $k$ (in case it is positive). For any variety $X$ over $k$, let   
  $$
  H^j(X,\QQ_l(i))
  $$ 
be the $j$-th $l$-adic \'etale cohomology group of the variety $X$ twisted by $i$. Then $H^*(-,\QQ _l)$ is a Weil cohomology theory over $k$. 

In particular, for any smooth projective $X$ over $k$, there is a cycle class homomorphism 
  $$
  \cl _X^j:CH^j(X)\to H^{2j}(X,\QQ_l(j))\; ,
  $$
whose kernel is the group 
  $$
  CH^j(X)_{\hom }
  $$
of algebraic cycles of codimension $j$ modulo homological equivalence on $X$.

If the ground field is $\CC $, then $H^j(X,\QQ_l(i))$ is isomorphic to the Betti cohomology group $H^j(X(\CC ),\QQ _l)$ with coefficients in $\QQ _l$, and embedding $\QQ _l$ in to $\CC $ (non-canonically) we obtain that the group $H^j(X(\CC ),\QQ _l)\otimes \CC $ is isomorphic to the group $H^j(X(\CC ),\CC )$, which can be useful if we want to compare \'etale to Hodge setting. 

Now, for any smooth projective connected variety $X$ of dimension $n$ over $k$ the class 
  $$
  cl ^n_{X\times X}(\Delta _X)\in H^{2n}(X\times X,\QQ _l)
  $$ 
of the diagonal $\Delta _X$ decomposes into its K\"unneth components 
  $$
  {cl(\Delta _X)}_{i,n-i}\in H^i(X,\QQ _l)\otimes H^{n-i}(X,\QQ _l)\; ,
  $$
for all 
  $$
  0\leq i\leq 2n\; .
  $$
It is a part of the Standard Conjectures on algebraic cycles that these cohomological classes can be lifted to mutually orthogonal idempotents 
  $$
  \pi _i\in CH^0(X,X)_{\QQ }\; ,
  $$
such that
  $$
  \Delta _X=\sum _{i=1}^{2n}\pi _i
  $$
in the associative algebra $CH^0(X,X)_{\QQ }$. 

In \cite{MurreIndag} Murre conjectured that the correspondences 
  $$
  \pi _0,\ldots ,\pi _{j-1}
  $$ 
and 
  $$
  \pi _{2j+1},\ldots ,\pi _{2n}
  $$ 
act as zero on $CH^j(X)_{\QQ }$, for any $0\leq j\leq n$, that the decreasing filtration
  $$
  F^iCH^j(X)_{\QQ }=
  \ker ({\pi _{2j}}_*)\cap \ker ({\pi _{2j-1}}_*)
  \cap \ldots \cap \ker ({\pi _{2j-i+1}}_*)
  $$
\noindent is independent on the choice of the projectors $\pi _0,\ldots ,\pi _{2n}$, and
  $$
  F^1CH^j(X)_{\QQ }=CH^j(X)_{\hom ,\QQ }
  $$

\medskip

\noindent for each $0\leq j\leq n$.

Murre's conjectures are equivalent to the conjectures of Beilinson and Bloch, taken for all smooth and projective $X$ over $k$, see \cite{JannsenMotivicSheaves}. 

For short, write
  $$
  M^i(X)=(X,\pi _i,0)\; ,
  $$
so that 
  $$
  M(X)=\oplus _{i=0}^{2n}M(X)
  $$
is the Chow-K\"unneth decomposition of the motive $M(X)$. 

Choose and fix a closed point $P_0$ on $X$. Then 
  $$
  \pi _0=[P_0\times X]\; ,
  $$
  $$
  \pi _{2n}=[X\times P_0]
  $$
and
  $$
  M^0(X)=\uno \; ,
  $$
  $$
  M^{2n}(X)=\Le ^n
  $$
in $\CHM (k,\QQ )$, and, in fact, in $\CHM (k,\ZZ )$. 

Let, for example, $C$ be a smooth projective curve over $k$. Then 
  $$
  \pi _1=\Delta _C-\pi _0-\pi _2
  $$
is also an integral projector, and we obtain the well-known decomposition
  \begin{equation}
  \label{bereza}
  M(C)=\uno \oplus M^1(C)\oplus \Le
  \end{equation}
in the category $\CHM (k,\ZZ )$. Murre's conjectures are true for curves. The motives $\uno $ and $\Le $ are evenly $1$-dimensional, and the motive $M^1(C)$ is oddly $2g$-dimensional in Kimura's sense, see \cite{Kimura}. Here $g$ is the genus of the curve $C$, 

If $n>1$, one can construct the Picard and its dual Albanese projector, $\pi _1$ and $\pi _{2n-1}$, which determine the Picard motive $M^1(X)$ and the Albanese motive $M^{2n-1}(X)$ respectively, both with coefficients in $\QQ $, which have the expected behaviour, see the details in \cite{MurreCrelle}.

Let $S$ be a smooth projective surface over $k$. Then
  $$
  \pi _2=\Delta _S-\pi _0-\pi _1-\pi _3-\pi _4
  $$
is the middle projector $\pi ^2$, and, respectively, we obtain the Chow-K\"unneth decomposition 
  \begin{equation}
  \label{olha}
  M(S)=\uno \oplus M^1(S)\oplus M^2(S)\oplus M^3(S)\oplus \Le ^2  \end{equation}
of the motive $M(S)$ in $\CHM (k,\QQ )$. 

This decomposition can be refined further by splitting the algebraic part from $M^2(S)$, see \cite{KMP}. 

Namely, let $\rho $ be the Picard number of $S$ and choose $\rho $ divisors 
  $$
  D_1,\ldots ,D_{\rho }
  $$ 
whose cohomology classes generate the second Weil cohomology group $H^2(S,\QQ _l)$. Choose the Poincar\'e dual divisors 
  $$
  D'_1,\ldots ,D'_{\rho }\; ,
  $$
 with coefficients in $\QQ $, so that the intersection number 
   $$
   \langle D_i\cdot D'_j\rangle 
   $$ 
is the Kronecker symbol. Then
  $$
  \pi ^{\alg }_2=\sum _{i=1}^{\rho }D_i\times D_i'
  $$
is a projector cutting the algebraic part in $M^2(S)$. That is, if 
  $$
  \pi ^{\tr }_2=\pi _2-\pi ^{\alg }_2\; ,
  $$
then we obtain the decomposition
  \begin{equation}
  \label{yasen}
  M^2(S)=M^2_{\alg }(S)\oplus M^2_{\tr }(S)
  \end{equation}
in $\CHM (k,\QQ )$. 

Notice also that if 
  $$
  H^1(S,\QQ _l)=H^3(S,\QQ _l)=0\; ,
  $$
i.e. if the surface $S$ is regular, then 
  $$
  M^1(S)=M^3(S)=0\; ,
  $$
and we obtain the decomposition 
  $$
  M(S)=\uno \oplus M^2(S)\oplus \Le ^2
  $$
already in the integral category $\CHM (k,\ZZ )$. However, the decomposition (\ref{yasen}) is in $\CHM (k,\QQ )$. 

The Murre conjectures are known to be true for surfaces, except for independence of the filtration on the choice of the projectors $\pi _i$, and the latter is true if the motive $M(S)$ is finite-dimensional. 

The Chow-K\"unneth decomposition of the product of two smooth projective varieties $Y_1$ and $Y_2$ can be constructed as the K\"unneth product of the Chow-K\"unneth decompositions for both factors, if they exists, see \cite{MurreIndag}. In particular, if 
  $$
  S=C_1\times C_2
  $$
is the product of two smooth projective curves over $k$, then then surface $S$ admits an integral Chow-K\"unneth decomposition in the category $\CHM (k,\ZZ )$ by the formula
  $$
  M^i(S)=\bigoplus _{p+q=i}(M^p(C_1)\otimes M^q(C_2))\; .
  $$

In dimension $3$ some partial results are obtained too. In \cite{MurreIndag} Murre studied the case when 
  $$
  X=S\times C\; ,
  $$
where $S$ is a smooth projective surface and $C$ is a smooth projective curve over $k$, see \cite{MurreIndag}. Naturally, if $H^1(S,\QQ _l)$ and $H^3(S,\QQ _l)$ vanish, then the Chow-K\"unneth decomposition is integral as the product of two integral decompositions.

The motive of a smooth projective Fano threefold is finite-dimensional and the explicit Chow-K\"unneth decomposition of such a motive is studied in \cite{GG-Fano}.

\subsection{Chow motives of cubic hypersurfaces in $\PR ^5$}

Let now $X$ be a smooth hypersurface in $\PR ^{n+1}$, and let 
  $$
  b_j(X)=\dim H^j(X,\QQ _l)
  $$
be the $j$-th Betti number of the variety $X$. Then $b_j(X)$ is $0$ if $j$ is odd and $j\neq n$, and $b_j(X)$ is $1$ if $j$ is even and $j\neq n$. Then all cohomology groups $H^{2j}(X)$ are algebraic, for $2j\neq n$. 

Let $Y$ be a general hyperplane section of $X$, and let 
  $$
  \gamma =[Y]
  $$ 
 be its class in $CH^1(X)$. For any number $j$ between $0$ and $n$ let 
  $$
  \gamma ^j=[Y]\cdot \ldots \cdot [Y]
  $$ 
be the $j$-fold self-intersection of the class $\gamma $ in $CH^j(X)$. By the Lefschetz hyperplane section theorem, the vector space $H^{2j}(X)$ is generated by the cycle class $\gamma ^j$, if $2j\neq n$. 

Next, for any integer 
  $$
  0\leq i\leq 2n
  $$ 
let

\medskip 

  $$
  \pi _i=
  \left\{
  \begin{array}{ll}
  0 & \mbox{if $i=2j+1$, $0\leq j\leq n-1$ and $i\neq n$} \\
  \frac{1}{\deg (X)}\cdot \gamma ^{n-j}\times \gamma ^j & \mbox{if $i=2j$, $0\leq j\leq n$ and $i\neq n$}
  \end{array}
  \right.
  $$

\medskip

\noindent and let

\medskip

  $$
  \pi _n=
  \Delta _X-\sum _{\substack {i=0 \\ i\neq n}}^{2n}\pi _i\; .
  $$

\medskip

\noindent Such defined correspondences $\pi _0,\ldots ,\pi _{2n}$ give us the Chow-K\"unneth decomposition of the diagonal for $X$. 

In particular, let $X$ be a smooth cubic fourfold hypersurface in $\PR ^5$ over $k$, and assume that $k$ is of characteristic $0$. 
 
 Since 
   $$
   \deg (X)<5\; ,
   $$
the hypersurface $X$ is rationally connected, whence
  $$
  CH_0(X)_{\QQ }=\QQ \; .
  $$

Fix a point $P_0$ on $X$. Then, as above,
  $$
  \pi _0=[P_0\times X]\; ,
  $$
  $$
  \pi _1=0\; ,
  $$
  $$
  \pi _2=\frac{1}{3}\cdot \gamma ^3\times \gamma \; ,
  $$
  $$
  \pi _3=0\; ,
  $$
  $$
  \pi _4=\Delta _X-
  \sum _{\substack {i=0 \\ i\neq 4}}^8\pi _i\quad \hbox{\footnotesize{(no explicite construction)}}\; ,
  $$
  $$
  \pi _5=0\; ,
  $$
  $$
  \pi _6=\frac{1}{3}\cdot \gamma \times \gamma ^3\; ,
  $$
  $$
  \pi _7=0
  $$
and
  $$
  \pi _8=[X\times P_0]\; .
  $$
This gives the Chow-K\"unneth splitting
  $$
  M(X)=
  \uno \oplus \Le ^2\oplus M^4(X)\oplus \Le ^6\oplus \Le ^8
  $$
in the rational category $\CHM (k,\QQ )$.

Let $\rho _2$ be the rank of the algebraic part in $H^4(X,\QQ _l)$. Choosing $2$-cycles
  $$
  D_1,\ldots ,D_{\rho _2}
  $$
on $X$, and their Poincar\'e dual $2$-cycles
  $$
  D'_1,\ldots ,D'_{\rho _2}\; ,
  $$
exactly in the same way as we did it for surfaces above, one can easily construct the splitting
  $$
  M^4(X)=M^4_{\alg }(X)\oplus M^4_{\tr }(X)\; ,
  $$
in $\CHM (k,\QQ )$, where
  $$
  M^4_{\alg }(X)=\Le ^{\oplus \rho _2} \; ,
  $$
i.e.
  $$
  \pi ^4_{\alg }=\sum _{i=1}^{\rho _2}[D_i\times D'_i]\; .
  $$

Clearly, each copy of the Lefschetz motive $\Le $ is the motive $(X,D_i\times D'_i,0)$, and the transcendental motive 
  $$
  M^4_{\tr }(X)
  $$ 
is given by the projector
  $$
  \pi ^4_{\tr }=\pi _4-\pi ^4_{\alg }\; .
  $$

Let also
  $$
  \pi ^4_{\prim } =\Delta _X-\frac{1}{3}\cdot
  \sum _{j=0}^4\gamma ^{4-j}\times \gamma ^j\; ,
  $$
and let
  $$
  M^4_{\prim }(X)=(X,\pi ^4_{\prim },0)
  $$
be the primitive part of the motive $M(X)$, see \cite{KimuraMurreConf}. If the cubic $X\subset \PR ^5$ is very general, the results in \cite{ZarhinBollettino} show that
  $$
  \rho _2=1\; ,
  $$
whence
  $$
  M^4_{\prim }(X)=M^4_{\tr }(X)\; .
  $$
Then, for a very general cubic $X$, we get
  $$
  M^4(X)=
  \Le ^{\oplus \rho _2}\oplus M^4_{\prim }(X)\; .
  $$

Moreover, if $X$ is very general, then
  $$
  \End _{\QQ }(H^4(X)_{\prim })=\QQ \; ,
  $$
i.e. the rational Hodge structure on the middle primitive cohomology is indecomposable, see Remark 2.6(a) in \cite{ZarhinReine} and Lemma 5.1 in \cite{VoisinUnivCH0}.

\section{Proofs of main theorems}

\subsection{Motivic version of Kulikov's reduction (Theorem A)}

Let again $X$ be a general cubic hypersurface in $\PR ^5$ over $\CC $. The aim of this section is to prove Theorem A in Introduction. Namely, we want to show that if we assume that a very general $X$ is rational, that would imply that the integral motive $M(S)$ of any sufficiently general smooth projective surface $S$ over $\CC $ is integrally indecomposable. 

\begin{theorem}
\label{TheoremA}
If the motivic indecomposability conjecture is true, and if the motive of any smooth projective surface is finite-dimensional, then a very general cubic fourfold hypersurface in $\PR ^5$ is not rational.
\end{theorem}

\begin{pf}
Suppose $X$ is rational, and consider the corresponding birational map
  $$
  \PR ^4\dashrightarrow X\; .
  $$
Resolving the indeterminacy locus, we get a regular dominant morphism
  $$
  f:Y\to X
  $$
over $k$, where $Y$ is obtained by a chain of blow up operations at points, curves and surfaces, starting from $\PR ^4$.

A crucial geometric argument is this. Let
  $$
  F=F(X)
  $$
be the Fano variety of the cubic $X$. By the result of Voisin, there exists a surface
  $$
  F_0\subset F\; ,
  $$
such that any two points on $F_0$ are rationally equivalent on the fourfold $F$, see \cite{Voisin}. Moreover, for any line $L$ on $X$, such that its class $[L]$ in $F$ sits on the surface $F_0$, the triple line $3L$ is rationally equivalent to the third intersection power,
  $$
  [3L]=\gamma ^3\; ,
  $$
of the general hyperplane section $\gamma $ of the cubic $X$, see Lemma A.3(v) in \cite{ShenVial}. 

It follows that the class $\gamma $ of the hyperplane section in $CH^1(X)$ is divisible by $3$. Therefore, the splitting
  $$
  M(X)=
  \uno \oplus \Le ^2\oplus M^4(X)\oplus \Le ^6\oplus \Le ^8
  $$
is integral, i.e. it already happens in the integral category $\CHM (k,\ZZ )$. 

The morphism $f$ is generically $1:1$ and dominant. Therefore, the composition 
  $$
  \Gamma _f\circ \Gamma _f^{\trp }
  $$ 
is the identity automorphism of $M(X)$, again, in the integral category $\CHM (k,\ZZ )$. 

In other words, $f$ yields the embedding
  $$
  f^*=\Gamma _f^{\trp }:M(X)\to M(Y)\; ,
  $$
which integrally splits $M(X)$ from $M(Y)$, and therefore
  $$
  M(Y)=f^*(M(X))\oplus N
  $$
in $\CHM (k,\ZZ )$, where $f^*(M(X))$ is the submotive in $M(Y)$ cut out by the projector $\Gamma _f^{\trp }\circ \Gamma _f$ on $Y$.

Suppose we sequentially blow up $s_0$ points, $s_1$ curves $C_1,\ldots ,C_{s_1}$ and $s_2$ surfaces $S_1,\ldots ,S_{s_2}$ over $k$. Then the latter motive splits integrally as
  $$
  M(Y)=M(\PR ^4)\oplus M_0\oplus M_1\oplus M_2\; ,
  $$
where
  $$
  M_0=
  \oplus _{i=1}^{s_0}
  (\Le \oplus \Le ^2\oplus \Le ^3)\; ,
  $$
\medskip
  $$
  M_1=
  (\oplus _{i=1}^{s_1}
  M(C_i))\otimes (\Le \oplus \Le ^2)
  $$
and
  $$
  M_2=\oplus _{i=1}^{s_2}M(S_i)\otimes \Le \; .
  $$

\bigskip

As it was shown in \cite{Kulikov}, there exists an index 
  $$
  i_0\in \{ 1,\ldots ,s_2\} \; ,
  $$
such that the pullback under the morphism $f$ of the transcendental Hodge structure of the cubic $X$, being twisted by $1$, is an integral sub-Hodge structure in the transcendental Hodge structure of $S_{i_0}$. More importantly, this integral sub-Hodge structure does not equal to the whole transcendental Hodge structure of $S_{i_0}$.

Next, consider the monoidal functor 
  $$
  \CHM (k,\ZZ )\to \Num (k,\ZZ )\; ,
  $$
introduced in Section \ref{basics}. Recall that it sends a Chow motive $M$ to the motive $\bar M$ given by the same projector but modulo numerical equivalence, and similarly on morphisms.

The integral splitting
  \begin{equation}
  \label{okolica}
  \bar M(Y)=f^*(\bar M(X))\oplus \bar N
  \end{equation}
in $\CHM (k,\ZZ )$ induces the integral splitting
  \begin{equation}
  \label{zhautki}
  \bar M^2_{\tr }(S_{i_0})=
  (f^*(\bar M^4_{\prim }(X))\otimes \Ta )\oplus
  (\bar N_{i_0}\otimes \Ta )
  \end{equation}
in the category $\Num (k,\QQ )$.

Now, the motives of curves are finite-dimensional by Theorem 4.2 in \cite{Kimura}. Assuming that motives of all smooth projective surfaces are finite-dimensional, we obtain, in particular, that the motives $M(S_i)$ are all finite-dimensional in Kimura's sense. Then the motive $M(Y)$ is finite-dimensional, and, of course, the motive $M(X)$ is also finite-dimensional, as objects in $\CHM (k,\QQ )$. 

As the cubic $X$ is very general in $\PR ^5$,
  $$
  \bar M^4_{\tr }(X)=\bar M^4_{\prim }(X)\; ,
  $$
and this motive is indecomposable by Lemma 5.1 in \cite{VoisinUnivCH0} and the absence of phantom submotives in finite-dimensional motives, which is due to Kimura's Proposition 7.5 in \cite{Kimura}. 

Lemma 3 in \cite{Kulikov} then gives us that
  $$
  \bar N_{i_0}\neq 0\; ,
  $$
so that both summands in (\ref{zhautki}) are nontrivial.

Let
  $$
  N^2(S_{i_0}\times S_{i_0})
  $$
be the group of $2$ cycles modulo numerical equivalence on $S_{i_0}\times S_{i_0}$. In terms of correspondences modulo numerical equivalence, the splitting (\ref{okolica}) induces a decomposition
  $$
  \bar \Delta =\bar \Lambda +\bar \Xi
  $$
of the diagonal class $\bar \Delta $ of the cubic $X$ in to two orthogonal idempotents in $N^2(S_{i_0}\times S_{i_0})$, such that
  $$
  f^*(\bar M^4_{\prim }(X))\otimes \Ta =M_{\bar \Lambda }
  \qquad \hbox{and}\qquad
  \bar N_{i_0}\otimes \Ta =M_{\bar \Xi }
  $$
in $\Num (k,\QQ )$.

Since the motive $M(S_{i_0})$ is finite-dimensional, all numerically trivial endomorphisms of $M(S_{i_0})$ are nilpotent by Proposition 7.5 in \cite{Kimura}. 

Then the standard lifting idempotent property gives us that there exist two orthogonal idempotents
  $$
  \Lambda '\; ,\; \Xi '\in CH^2(S_{i_0}\times S_{i_0})\; ,
  $$
such that
  $$
  \bar \Lambda '=\bar \Lambda \; ,\quad
  \bar \Xi '=\bar \Xi
  $$
and
  $$
  \Delta =\Lambda +\Xi
  $$
in $CH^2(S_{i_0}\times S_{i_0})$. Therefore, we may assume that $\Lambda $ and $\Xi $ are orthogonal idempotents from the very beginning, and we obtain the corresponding integral decomposition
  \begin{equation}
  \label{senoval}
  M^2(S_{i_0})=M_{\Lambda }\oplus M_{\Xi }
  \end{equation}
in $\CHM (k,\ZZ )$, such that
  $$
  \bar M_{\Lambda }=M_{\bar \Lambda }
  $$
and
  $$
  \bar M_{\Xi }=M_{\bar \Xi }\; .
  $$

The integral decomposition (\ref{senoval}) remains a non-trivial decomposition modulo numerical equivalence with coefficients in $\ZZ $, and both projectors $\Lambda $ and $\Xi $ are not balanced modulo rational equivalence with coefficients in $\QQ $.

As the threefold $X$ is defined over the function field $\CC (\PR ^{55})$, the minimal field of definition of the surface $S_{i_0}$ is of a big transcendence degree. Then we get a contradiction with the indecomposability assumption.
\end{pf}

\subsection{Lifting theorem over a Henselian DVR (Theorem B)}

Let us also prove Theorem B in Introduction. 

Let $A$ be a Henselian DVR with maximal ideal $\gm $, the ring of fractions
  $$
  K=A_{(0)}
  $$
and the residue field 
  $$
  \kappa =A/\gm \; .
  $$
Let 
  $$
  S=\Spec (A)\; ,
  $$
let 
  $$
  0=\Spec (\kappa )\qqand \eta =\Spec (K)
  $$
be the closed and generic point of $S$ respectively. 

Let
  $$
  X\to S
  $$
be a smooth projective family of surfaces over $S$ with the generic fibre 
  $$
  X_{\eta }=X_K
  $$ 
and the closed fibre 
  $$
  X_0=X_{\kappa }\; .
  $$ 
  
Let also $\bar \kappa $ and $\bar K$ be the algebraic closures of the fields $\kappa $ and $K$ respectively. The goal of this section is to prove the following

\begin{theorem}
\label{TheoremB}
Suppose that the motive $M(X_{\bar K})$ is finite-dimensional. Then, if the motive $M(X_{\bar \kappa })$ is integrally essentially indecomposable, the motive $M(X_{\bar K})$ is also integrally essentially indecomposable in $\CHM (k,\ZZ )$.
\end{theorem}

\begin{pf}
The generic and closed fibres of the smooth projective family 
  $$
  X\times _SX\to S
  $$
are the products $X_{\eta }\times _{\eta }X_{\eta }$ and $X_{\kappa }\times X_{\kappa }$ respectively, and the same for the geometric generic and closed fibres. Since $A$ is Henselian, we have the specialization commutative square


  $$
  \diagram
  CH^2(X_{\bar K}\times _{\bar K}X_{\bar K}) \ar[dd]_-{\spc } 
  \ar[rr]^-{\cl } & & H^4(X_{\bar K}\times _{\bar K}X_{\bar K},\QQ _l(2)) \ar[dd]^-{\spc }_{\simeq } \\ \\
  CH^2(X_{\bar \kappa }\times _{\bar \kappa }X_{\bar \kappa }) 
  \ar[rr]^-{\cl } & & H^4(X_{\bar \kappa }\times _{\bar \kappa }X_{\bar \kappa }),\QQ _l(2))
  \enddiagram
  $$

\bigskip

\noindent on Chow groups and cohomology, see Expose X, 7.13 (pp 576 - 581) in \cite{SGA6} and Example 20.3.5 on page 400 in \cite{FultonRatEqSingVar}. The horizontal homomorphisms in this diagram are cycle class homomorphisms to $l$-adic cohomology groups, and we assume that $l$ is coprime to the characteristics of both fields, $K$ and $\kappa $. 

The vertical arrows are the specialization homomorphisms, see Section 20.3 in Fulton's book \cite{Fulton}. Moreover, the right specialization is an isomorphism on cohomology groups, see Expose X, 7.13 in \cite{SGA6}. 

Assume now that the motive $M(X_{\bar \kappa })$ is integrally essentially indecomposable over $\bar \kappa $, but the motive $M(X_{\bar K})$ is integrally essentially decomposable over $\bar K$. Let
  $$
  \Delta _{X_{\bar K}}=\Lambda _{\bar K}+\Xi _{\bar K}
  $$
be the corresponding decomposition of the diagonal in $CH^2(X_{\bar \eta }\times X_{\bar \eta })$ in to the sum of two integral mutually orthogonal essential projectors on the surface $X_{\bar \eta }$, such that both projectors $\Lambda _{\bar K}$ and $\Xi _{\bar K}$ are homologically non-trivial and not balanced modulo rational equivalence with coefficients in $\ZZ $ on $X_{\bar K}\times _{\bar K}X_{\bar K}$. 

The specialization of the diagonal class $\Delta _{\bar K}$ is the diagonal class $\Delta _{\bar \kappa }$. Let
  $$
  \Delta _{X_{\bar \kappa }}=\Lambda _{\bar \kappa }+\Xi _{\bar \kappa }
  $$
be the image of the decomposition of the diagonal $\Delta _{\bar K}$ under the left specialization on Chow groups. 

Since the square above is commutative, and the right specialization is an isomorphism, both projectors $\Lambda _{\bar \kappa }$ and $\Xi _{\bar \kappa }$ are cohomologically non-trivial on $X_{\bar \kappa }\times _{\bar \kappa }X_{\bar \kappa }$.

On the other hand, since we assume that the motive $M(X_{\bar \kappa })$ is integrally essentially indecomposable, it follows that at least one of the projectors $\Lambda _{\bar \kappa }$ or $\Xi _{\bar \kappa }$ is balanced modulo rational, and hence modulo homological equivalence on $X_{\bar \kappa }\times _{\bar \kappa }X_{\bar \kappa }$. 


Suppose, $\Xi _{\bar \kappa }$ is balanced modulo homological equivalence on $X_{\bar \kappa }\times _{\bar \kappa }X_{\bar \kappa }$. Since the coniveau filtration is preserved by specialization, see pp 283 - 284 in \cite{Gillet}, it follows that $\Xi _{\bar K}$ is balanced modulo cohomological equivalence on $X_{\bar K}\times _{\bar K}X_{\bar K}$.

In other words, there exist a homologically trivial $2$-cyckle class $\gA _{\bar K}$ and a balanced $2$-cycle class $\gB _{\bar K}$, such that
  $$
  \Xi _{\bar K}=\gA _{\bar K}+\gB _{\bar K}
  $$
in $CH^2(X_{\bar K}\times _{\bar K}X_{\bar K})$. As we assume that the motive $M(X_{\bar K})$ is finite-dimensional, the correspondence $\gA _{\bar K}$ is nilpotent in the ring $CH^2(X_{\bar K}\times _{\bar K}X_{\bar K})_{\QQ }$. Since balanced correspondences is an ideal in this ring, we obtain that for a positive integer $N\gg 0$,
  $$
  \Xi _{\bar K}=\Xi _{\bar K}^N=(\gA _{\bar K}+\gB _{\bar K})^N=
  \gA _{\bar K}^N+\hbox{balanced correspondence}\; ,
  $$
and if $N$ is sufficiently big, we obtain that $\Xi _{\bar K}$ is balanced modulo rational equivalence with coefficients in $\QQ $ on $X_{\bar K}\times _{\bar K}X_{\bar K}$. This is a contradiction, and hence the motive $M(X_{\bar K})$ is essentially indecomposable. 
\end{pf}

\subsection{The Dirichlet theorem and unirationality of surfaces in $\PR ^3_{\bar \FF _p}$}

Recall the Dirichlet theorem on prime numbers in an arithmetical progression:

\begin{theorem}
\label{Dirichlet}
Let $a$ and $m$ be two integers, and assume they are coprime. Then there exists an infinite number of primes $p$, such that
  $$
  p\equiv a(\md m)\; .
  $$
\end{theorem}

The proof of this famous result can be found, for example, in \cite{IrelandRosen}.

\medskip
In particular, for any positive integer $n$, there exists an infinite number of primes $p$, such that
  $$
  p\equiv n-1(\md n)\; ,
  $$
and therefore
  $$
  p\equiv -1(\md n)\; .
  $$
By Shioda's main theorem in \cite{Shioda74}, the Fermat type hypersurface
  $$
  F_n:x_0^n+x_1^n+x_2^n+x_3^n=0
  $$
is unirational in $\PR ^3$, if the characteristic of the ground field is $p$. Thus, there exists an infinite number of primes $p$, such that the reduction of $F_n$ mod $p$ is good and unirational. 

A natural expectation here is that the integral motive of such a hypersurface is essentially indecomposable in $\CHM (\bar \FF _p,\ZZ )$. If this is true, then we can use Theorem \ref{TheoremA} to lift integral essential indecomposability to zero characteristic, and then lift it further to general hypersurfaces over the filed of functions on the parameter space, in order to prove motivic indecomposability conjecture stated in Introduction.

\bigskip

\begin{small}

\end{small}

\bigskip

\bigskip

\begin{small}

{\sc Department of Mathematical Sciences, University of Liverpool, Peach Street, Liverpool L69 7ZL, England, UK}

\end{small}

\medskip

\begin{footnotesize}

{\it E-mail address}: {\tt vladimir.guletskii@liverpool.ac.uk}

\end{footnotesize}

\end{document}